\documentclass[12pt]{amsart}

\usepackage{array,float, multirow}
\usepackage{amssymb, amsmath, stmaryrd,amsthm}

\theoremstyle{plain}
\newtheorem{theorem}{Theorem}[section]
\newtheorem{lemma}[theorem]{Lemma}
\newtheorem{corollary}[theorem]{Corollary}
\newtheorem{proposition}[theorem]{Proposition}

\theoremstyle{definition}

\theoremstyle{remark}

\newtheorem*{remark}{Remark}
\newtheorem*{acknowledgements}{Acknowledgements}
\newtheorem*{organisation/notation}{Organisation and notation}

\long\def\symbolfootnote[#1]#2{\begingroup%
\def\thefootnote{\fnsymbol{footnote}}\footnote[#1]{#2}\endgroup}

\def \Proj{\mathbb{P}}
\def \N {\ensuremath{\mathbb{N}}}

\def \Qp {\mathbb{Q}_p}
\def \Z {\mathbb{Z}}

\def \Zp  {\mathbb{Z}_p}
\def \Zpn {\mathbb{Z}_p^n}

\def \Fp  {\mathbb{F}_p}

\def \bfi {{\bf i}}
\def \bfj {{\bf j}}
\def \bfk {{\bf k}}

\def \bfr {{\bf r}}

\def \bfx {{\bf x}}

\def \mcN {{N}}

\def \GL {\text{GL}}

\def \tud {\textup{d}}
\def \tuL {\textup{L}}
\def \tut {\textup{t}}
\def \mcR {\mathcal{R}}
\def \diag {\textup{diag}}
\def \bfy {{\bf y}}
\def \G {\Gamma}
\def \ni {\noindent}
\renewcommand{\epsilon}{\varepsilon}
\renewcommand{\phi}{\varphi}

\begin{document}

\title[Zeta functions of $3$-dimensional $p$-adic Lie algebras]{Zeta
functions of $3$-dimensional\\ $p$-adic Lie algebras}

\author{Benjamin Klopsch}
\author{Christopher Voll}

\address{Benjamin Klopsch, Department of Mathematics, Royal Holloway,
University of London, Egham TW20 0EX, United Kingdom}
\email{Benjamin.Klopsch@rhul.ac.uk} 
\address{Christopher Voll, School of Mathematics, University of Southampton, University Road, Southampton SO17 1BJ, United Kingdom}
\email{C.Voll.98@cantab.net} \thanks{\today}


\begin{abstract} We give an explicit formula for the subalgebra zeta
  function of a general $3$-dimensional Lie algebra over the $p$-adic
  integers~$\Zp$. To this end, we associate to such a Lie algebra a
  ternary quadratic form over~$\Zp$. The formula for the zeta function
  is given in terms of Igusa's local zeta function associated to this
  form.
\end{abstract}

\maketitle

%-----------------------------------------------------------------------
\section{Introduction}
For a prime $p$, let $\Zp$ denote the ring of $p$-adic integers. The
 \emph{(subalgebra) zeta function} of a $\Zp$-algebra~$L$, additively
 isomorphic to $\Zp^n$, is the Dirichlet series
\begin{equation*}%`\label{definition zeta function}
 \zeta_{L}(s)=\sum_{H\leq L}|L:H|^{-s},
\end{equation*}
where the sum ranges over the subalgebras of finite index in~$L$, and
$s$ is a complex variable.

Zeta functions of $\Zp$-\emph{Lie} algebras play an important role in
the subject of subgroup growth. Indeed, to every saturable $p$-adic
analytic pro-$p$ group $G$ there is an associated $\Zp$-Lie algebra
$L=L(G)$ and if $\dim(G)\leq p$ then
$$\zeta_{L(G)}(s)=\zeta_{G}(s):=\sum_{H\leq G}|G:H|^{-s},$$ the
\emph{(subgroup) zeta function} of the group~$G$
(cf.~\cite{Klopsch/05, KlopschGonzalez-Sanchez/07} and references
therein). Similarly, for a finitely generated nilpotent group $G$,
there is a nilpotent $\Z$-Lie algebra $L(G)$ such that, for almost all
primes $p$,
$$\zeta_{\Zp\otimes_\Z L(G)}(s)=\zeta_{G,p}(s):=\sum_{H\leq_p
G}|G:H|^{-s},$$ the \emph{local (subgroup) zeta function} of $G$ at
the prime $p$, enumerating subgroups of finite $p$-power index in~$G$
(cf.~\cite{GSS/88}). Thus, to some degree, the study of subgroup zeta
functions reduces to the study of subalgebra zeta functions of
$\Zp$-Lie algebras. Zeta functions of groups and rings in general have
attracted considerable interest over the last few decades; we refer
to~\cite{duSG/06} for a recent survey.

By now, numerous examples of zeta functions of nilpotent and soluble
$\Zp$-Lie algebras have been calculated. One of the first examples is
the zeta function of the Heisenberg Lie algebra $\mathfrak{h}(\Zp)$
which was computed in~\cite{GSS/88}. For further explicit calculations
see, for example, \cite{Taylor/01, VollBLMS/06,
duSautoyWoodward/06}. On the other hand, explicit examples of zeta
functions of insoluble $\Zp$-Lie algebras are thin on the ground. Only
for zeta functions associated to the two $\Qp$-forms of the simple Lie
algebra of type $A_1$ have explicit formulae been found: using results
of Ilani~(\cite{Ilani/99}), du Sautoy gave a formula for the zeta
function of $\mathfrak{sl}_2(\Zp)$ (\cite{duSsl2/00}; see
also~\cite{duSTaylor/02}). In \cite{Klopsch/03}, Klopsch computed the
zeta function of $\mathfrak{sl}_1(\Delta_p)$, where $\Delta_p$ denotes
the maximal $\Zp$-order in a central $\Qp$-division algebra of
index~$2$. No explicit formula for the zeta function of any
`semi-simple' $\Zp$-Lie algebra of dimension greater than $3$ is known
(cf. \cite[Problem 9(c) on p. 431]{LubotzkySegal/03}).

In~\cite{Voll/06a}, Voll introduced a method for computing zeta
functions of $\Zp$-algebras in terms of certain $p$-adic integrals
generalising Igusa's local zeta function. Given a polynomial
$f(\bfx)\in\Zp[x_1,\dots,x_n]$, Igusa's local zeta function associated
to $f$ is defined as the $p$-adic integral
\begin{equation*}
 Z_f(s)=\int_{\Zp^n}|f(\bfx)|_p^s\tud\mu.
\end{equation*}
Here $|\;|_p$ denotes the $p$-adic absolute value, $s$ is a complex
variable and $\tud\mu$ stands for the normalised additive Haar measure
on~$\Zpn$. Igusa's local zeta function is closely connected to the
Poincar\'e series enumerating the numbers of solutions of the
congruences $f(\bfx)\equiv 0\mod(p^m)$ for $m\in\N$ (see
Section~\ref{igusa} for further details and references).

The purpose of the current paper is to demonstrate that this point of
view may be used to unify and generalise the existing computations of
zeta functions of $3$-dimensional $\Zp$-Lie algebras. These specific
calculations draw on a variety of methods: they range from elementary
counting arguments for the Heisenberg Lie algebra in~\cite{GSS/88}
over a carefully chosen resolution of singularities for a
high-dimensional hypersurface in \cite{duSTaylor/02} to a structural
analysis of a division algebra in~\cite{Klopsch/03}. Our main result
generalises these results and subsumes them under a unified
description in terms of a rather tame $p$-adic integral: Igusa's local
zeta function of a ternary quadratic form, naturally associated to a
$3$-dimensional $\Zp$-Lie algebra. To formulate our result we recall
that the zeta function of the abelian $3$-dimensional $\Zp$-Lie
algebra is
$$\zeta_{\Zp^3}(s)=\zeta_p(s)\zeta_p(s-1)\zeta_p(s-2).$$ Here
$\zeta_p(s)=(1-p^{-s})^{-1}$ denotes the $p$-th local Riemann zeta
function.

\begin{theorem}\label{theorem}
  Let $L$ be a $3$-dimensional $\Zp$-Lie algebra. Then there is a
  ternary quadratic form $f(\bfx)\in\Zp[x_1,x_2,x_3]$, unique up to
  equivalence, such that, for $i\geq 0$,
  \begin{equation*}%\label{zeta formula}
   \zeta_{p^iL}(s)=\zeta_{\Zp^3}(s)-Z_f(s-2)\zeta_p(2s-2)\zeta_p(s-2)p^{(2-s)(i+1)}(1-p^{-1})^{-1},
  \end{equation*}
 where $Z_f(s)$ is Igusa's local zeta function associated to $f$.
\end{theorem}

In the course of the proof of Theorem~\ref{theorem} we define
$f(\bfx)$ in terms of the structure constants of $L$ with respect to a
given $\Zp$-basis; different bases give rise to equivalent quadratic
forms (see Section~\ref{proof} for details).

The following table lists the ternary quadratic forms controlling the
subalgebra growth in several special cases mentioned above.

\begin{figure}[H]
 \begin{center}
  \begin{tabular}{|c|c|}\hline
  \textup{Lie algebra} & \textup{ternary quadratic form $f(\bfx)$}
  \\\hline $\Zp^3$&0\\ $\mathfrak{h}(\Zp)$&$x_3^2$\\
  $\mathfrak{sl}_2(\Zp)$&$x_3^2+4x_1x_2$\\ $\mathfrak{sl}_1(\Delta_p)$
  &$\begin{cases}-2(x_1^2-2(x_2^2-x_2x_3+x_3^2))&\text{ for }p=2\\
  x_3^2-\rho x_2^2-px_1^2 \text{ ($\rho$ a non-square mod $p$)}&\text{
  for $p$ odd}\end{cases}$\\\hline
  \end{tabular}
 \end{center}
\end{figure}

It is comparatively easy to compute Igusa's local zeta functions
associated to these forms. By Theorem \ref{theorem} we immediately
recover the known formulae for the zeta functions of the Lie algebras
in this table; see Sections~\ref{heisenberg}, \ref{sl2} and
\ref{D_p}. In Section \ref{section soluble} we use Theorem
\ref{theorem} to treat the soluble case, which has not been previously
studied. We give a complete list of the binary quadratic forms for the
soluble $3$-dimensional $\Zp$-Lie algebras (for odd $p$) and derive
formulae for Igusa's local zeta functions associated to these
quadratic forms. This leads to formulae for the zeta functions of all
these Lie algebras. Our computations show, in particular, that many
among the soluble $\Zp$-Lie algebras are isospectral
(i.e. non-isomorphic but sharing the same zeta function).

As an immediate corollary to Theorem~\ref{theorem} we gain control
over the real parts of the poles of zeta functions of $3$-dimensional
$\Zp$-Lie algebras. Indeed, it is easy to see that the only candidate
poles of Igusa's local zeta function of a ternary quadratic form have
real part $-3/2$, $-1$ and $-1/2$ (see Lemma~\ref{lemma IGLZ
  poles}). Thus we obtain

\begin{corollary}\label{coro 1}
 If $s$ is a pole of $\zeta_{p^iL}(s)$, then
$\mathfrak{Re}(s)\in\{0,1/2,1,3/2,2\}$.
\end{corollary}

\noindent To determine the poles of zeta functions of $\Zp$-(Lie)
algebras in general is a difficult and almost entirely unsolved
problem. Of particular importance is the largest actually occurring
real pole, as its position and order determine the asymptotics of the
subalgebra growth of~$L$. Our analysis allows us to solve this problem
for the soluble $\Zp$-Lie algebras of dimension~$3$ (see
Proposition~\ref{prop} and Corollary~\ref{coro}).

\medskip
\ni We conclude the introduction with a number of remarks.

\smallskip
\ni 1. Though our results are formulated for Lie algebras, our
arguments only draw on the fact that a Lie algebra is
\emph{antisymmetric}; the Jacobi identity is not being used anywhere.

\smallskip
\ni 2. In general, no simple identity is known which relates the zeta
function of a $d$-dimensional $\Zp$-algebra $L$ with that of $p^iL$,
$i\in\N_0$. In~\cite[Theorem 2.1]{duSsl2/00}, du~Sautoy gives such a formula for
the special case $d=3$; our Theorem~\ref{theorem} provides inter alia
an alternative proof of this formula, without reference to Mann's work
on probabilistic zeta functions~(\cite{Mann/96}).

\smallskip
\ni 3. We point out that Theorem~\ref{theorem} is a `local result',
whereas the main conclusions of the results in~\cite{Voll/06a} are
valid for almost all completions of a `global object' (such as a
torsion-free nilpotent group or a torsion-free ring). The main
application of the approach developed in~\cite{Voll/06a} is to prove
that, given a torsion-free ring~$L$ (not necessarily associative or
Lie), the associated local zeta function $\zeta_{\Zp\otimes_\Z L}(s)$
satisfies a functional equation for almost all primes $p$. The
occurrence of functional equations for the zeta functions of some
$3$-dimensional $\Zp$-Lie algebras is therefore only explained by the
results of~\cite{Voll/06a} in case these algebras are the `generic'
completions of an algebra over~$\Z$ (such as, for example,
$\mathfrak{sl}_2(\Zp)$ for odd~$p$). This corresponds to the fact that
the proof of a functional equation for Igusa's local zeta function
given in~\cite{DenefMeuser/91} critically depends on good reduction
modulo~$p$.

\smallskip \ni 4. The case of $3$-dimensional Lie algebras is the
first non-trivial one as far as subalgebra zeta functions are
concerned. In dimensions $1$ and $2$ it is not hard to see that the
concepts of subalgebra and additive sublattice coincide, so that the
zeta functions coincide with the well-known zeta functions for the
abelian case. The work in~\cite{Voll/06a} makes essential use of
generalisations of Igusa's local zeta functions to polynomial mappings
and several variables. It is remarkable that Theorem~\ref{theorem},
however, reduces the case of a $3$-dimensional Lie algebra to the
computation of the classical Igusa integral associated of a single
ternary quadratic form. Things get radically more complicated in
dimensions greater than $3$.

\smallskip \ni 5. Rather than counting all subalgebras of finite index
in a ring $L$, one may restrict attention to subalgebras with
additional algebraic properties; among the variants of $\zeta_L(s)$
that have been considered are the ideal zeta function
$\zeta^\triangleleft_L(s)$, enumerating ideals of finite index, and
the zeta function $\zeta^{\widehat{\;}}_L(s)$, counting subalgebras
isomorphic to~$L$. It would be interesting to study these zeta
functions of $3$-dimensional Lie algebras with the methods introduced
in the present paper.

\begin{organisation/notation}
 We prove our main result in
Section~\ref{proof}. In~Section~\ref{igusa} we collect a few
elementary observations about Igusa's local zeta function. The
examples given in the above table are studied in detail in
Section~\ref{examples}, where we derive the known formulae for their
zeta functions using Theorem~\ref{theorem}. In this section we also
list the binary quadratic forms associated to $3$-dimensional soluble
$\Zp$-algebras ($p\geq3$) and compute their zeta functions.

Throughout this paper we denote by $\N$ the set of positive integers
and by $\N_0$ the set of non-negative integers. Given a prime $p$, we
denote by $\Zp$ the ring of $p$-adic integers and by $\Qp$ the field
of $p$-adic numbers. We write $v_p$ for the $p$-adic valuation, and
$|\;|_p$ for the $p$-adic absolute value.
\end{organisation/notation}

\section{Preliminaries on Igusa's local zeta function} \label{igusa}

For general background on the theory of Igusa's local zeta function we
refer the reader to~\cite{Denef/91, Igusa/00}. Let
$f(\bfx)\in\Zp[x_1,\dots,x_n]$. The well-known connection between
Igusa's local zeta function $Z_f(s)$ and the Poincar\'e series
enumerating the numbers of solutions of $f(\bfx)\equiv 0$ modulo
$(p^m)$ mentioned in the introduction is the following. For
$m\in\N_0$, set
\begin{equation*}
 N_m:=|\{\bfx\in(\Z/p^m\Z)^n\mid f(\bfx)=0\}|.
\end{equation*} 
We write $t=p^{-s}$ and treat it as an independent variable. A simple
computation (cf.~\cite[Section 2.1]{Denef/91}) shows that the Poincar\'e series
\begin{equation*}
 P_f(t):=\sum_{m=0}^\infty N_m (p^{-n}t)^{m}
\end{equation*}
is related to Igusa's local zeta function by the formula
\begin{equation}\label{relation poincare igusa}
 P_f(t)=\frac{1-tZ_f(s)}{1-t}.
\end{equation}

Whilst the Poincar\'e series $P_f(t)$ counts solutions of polynomial
equations in \emph{affine} space we shall see in Section~\ref{proof}
below that counting subalgebras is related to counting solutions of
polynomial equations in finite \emph{projective} spaces. To this end
we define now, for $m\in\N_0$, the affine cones
\begin{align*}
 W&:=\Zpn\setminus p\Zpn,\\\quad W(m)&:=(W+(p^m\Zp)^n)/(p^m\Zp)^n
\end{align*}
and set $$N_m^\star:=|\{\bfx\in W(m)\mid
f(\bfx)=0\}|.$$ We will utilise a formula for
\begin{equation}\label{Pstar formula}
 P_f^\star(t):=\sum_{m=0}^\infty N_m^\star (p^{-n}t)^{m},
\end{equation}
analogous to~\eqref{relation poincare igusa} in the special case that
$f$ is \emph{homogeneous}. If $f$ is homogeneous of degree $d$, say,
we have (cf.~\cite[(1) on p. 1141]{DenefMeuser/91})
\begin{equation}\label{Z Zstar}
 Z_f(s)=\frac{1}{1-p^{-n-ds}}Z_f^\star(s),
\end{equation}
where 
\begin{equation*}
 Z_f^\star(s)=\int_{W}|f(\bfx)|_p^s\tud\mu.
\end{equation*}

\begin{lemma} \label{Pstar lemma}If $f$ is homogeneous, then
 \begin{equation}\label{Pstar Zstar}
  P_f^\star(t)=\frac{1-p^{-n}t-tZ_f^\star(s)}{1-t}.
 \end{equation}
\end{lemma}
\begin{proof}
 For $m\in\N_0$, set $\mu_m^\star:=\mu(\{\bfx\in W\mid
 v_p(f(\bfx))=m\}).$ We claim that
 \begin{equation}\label{mustar formula}
  \mu_m^\star=\frac{N^\star_m}{p^{nm}}-\frac{N^\star_{m+1}}{p^{n(m+1)}}-\delta_{m,0}\,p^{-n},
 \end{equation}
 where $\delta_{m,0}$ denotes the Kronecker-delta. Indeed, for $m\in\N_0$ we may write 
 \begin{equation*}
 \mu_m^\star=\mu(\{\bfx\in W\mid v_p(f(\bfx))\geq m\})-\mu(\{\bfx\in
 W\mid v_p(f(\bfx))\geq m+1\}).
 \end{equation*}
We have 
\begin{equation*}
\mu(\{\bfx\in W\mid v_p(f(\bfx))\geq
m\})= 
\begin{cases} N^\star_m/p^{nm} & \text{if $m\geq 1$,}
\\
\mu(W)=1-p^{-n} &\text{if $m=0$.}
\end{cases}
\end{equation*}
Using~\eqref{mustar formula} we obtain
 \begin{align}
  Z_f^\star(s)&=\sum_{m=0}^\infty \mu_m^\star p^{-ms}=\sum_{m=0}^\infty
            \left(\frac{N^\star_m}{p^{nm}}-\frac{N^\star_{m+1}}{p^{n(m+1)}}\right)p^{-ms}-p^{-n}\label{Zstar formula}\\
            &=P_f^\star(t)-\frac{1}{t}\left(P_f^\star(t)-1\right)-p^{-n}.\nonumber
 \end{align}
The lemma follows.
\end{proof}

To describe the position of the poles of Igusa's local zeta function
is in general a difficult and interesting problem. In the current
paper we shall only work with Igusa's local zeta function associated
to ternary quadratic forms. This case is well-understood:

\begin{lemma} \label{lemma IGLZ poles}
Let $f(\bfx)\in\Zp[x_1,x_2,x_3]$ be a ternary quadratic form. If $s$
is a pole of $Z_f(s)$, then $\mathfrak{Re}(s)\in\{-3/2,-1,-1/2\}$.
\end{lemma}

\begin{proof} 
  A ternary quadratic form defines a cone in affine $3$-space over a
  (possibly anisotropic) conic. A resolution of singularities is
  achieved by blowing up the origin, yielding an exceptional divisor
  with numerical data $(\nu,N)=(3,2)$. The divisors of the proper
  transform have numerical data $(1,1)$ unless the conic is a double
  line, in which case the numerical data is $(1,2)$. The real parts of
  the poles are to be found among the fractions $-\nu/N$
  (cf.~\cite{Denef/91, Igusa/75} for details).
\end{proof}
We thank Wim Veys for pointing this out to us.

\section{Proof of the main result}\label{proof}

In this section we prove Theorem~\ref{theorem}. Let $L$ be a
$3$-dimensional $\Zp$-Lie algebra. To compute the zeta function of $L$
it is helpful to make the following observations: the homothety
class $[\Lambda]$ of any $\Zp$-lattice $\Lambda$ in the $\Qp$-vector
space $\Qp\otimes_{\Zp}L$ contains a unique ($\subseteq$-) maximal
sub\emph{algebra} $\Lambda_0$ of~$L$, and the subalgebras contained in
the class $[\Lambda]$ are exactly the multiples $p^i\Lambda_0$,
$i\in\N_0$. Thus
\begin{equation*}%\label{reduction}
 \zeta_{L}(s)=\frac{1}{1-p^{-3s}}A(s), \text{ where
 }A(s):=\sum_{[\Lambda]}|L:\Lambda_0|^{-s}.
\end{equation*}
For the computation of $A(s)$ it is useful to sort the lattice classes
$[\Lambda]$ by their elementary divisor type with respect to the class
$[L]$, and to take advantage of the transitive action of the
group~$\G:=\textup{GL}_3(\Zp)$ on the classes of any fixed elementary
divisor type.

Write $L=\Zp e_1\oplus\Zp e_2\oplus\Zp e_3$. A sublattice
$\Lambda\subseteq L$ corresponds to a right-coset $\G M$, where
$M\in\GL_3(\Qp)\cap\textup{Mat}_3(\Zp)$, the set of
$3\times3$-matrices over $\Zp$ with non-zero determinant: the lattice
is generated by vectors whose coordinates with respect to the chosen
basis $(e_1,e_2,e_3)$ are encoded in the rows of~$M$. By the
elementary divisor theorem the right-coset $\G M$ contains a
representative of the form $D\alpha^{-1}$, where $\alpha\in\G$ and
$D=\diag(D_1,D_2,D_3)=p^{r_0}\diag(p^{r_1+r_2},p^{r_2},1)$ is a
diagonal matrix with $r_i\in\N_0$ for $i\in\{0,1,2\}$. We say that the
homothety class $[\Lambda]$ is of type\footnote{Note that this
  terminology differs slightly from the one used in \cite{Voll/06a},
  where the shape of diagonal matrix $D$ is encoded in a subset
  $I\subseteq[n-1]$ and a \emph{positive} vector $(r_i)\in\N^{|I|}$.}
$\bfr=(r_1,r_2)\in\N_0^2$ if the diagonal matrix $D$ determined by
$\Lambda$ is a scalar multiple of
$\diag(p^{r_1+r_2},p^{r_2},1)$. Below we shall make a case distinction
with respect to the invariant
$I([\Lambda])=\{i\in\{1,2\}|\,r_i\not=0\}$.

The matrix $\alpha$ is determined only up to right-multiplication by
an element of $\G_\bfr:=\textup{Stab}_{\G}(\G D)$, the stabiliser in
$\G$ of the right-coset $\G D$ under right-multiplication. The various
stabilisers will be described in detail below. A lattice class
$[\Lambda]$ is thus given by the pair $\bfr\in\N_0^2$ and a left-coset
$\alpha\G_\bfr\in\G/\G_\bfr$.

This parametrisation allows us to give a convenient description of the
index $|L:\Lambda_0|$ of the maximal subalgebra~$\Lambda_0$ in the
homothety class~$[\Lambda]$. In order to decide whether a lattice
$\Lambda$ is a subalgebra it suffices to check whether the products of
pairs of a given set of generators are contained in $\Lambda$. This is
particularly easy to verify if the right-coset $\G M$ contains a
diagonal matrix; in this case the condition of being a subalgebra
translates into a set of divisibility conditions on quadratic
polynomials in the entries of~$M$. In general, however, the coset~$\G
M$ may not contain any diagonal element. In this case a base change --
effectuated by right-multiplication with an element in the left-coset
$\alpha\Gamma_{\bfr}$ -- brings us into this desirable
situation. (Note that this approach differs distinctly from the point
of view taken e.g. in~\cite{duSG/00}, where all calculations are
performed with respect to a fixed basis.) More precisely, as indicated
in~\cite[Section 3]{Voll/06a}, a lattice $\Lambda$ corresponding to a
right-coset $\Gamma D \alpha^{-1}$ is a subalgebra of $L$ if and only
if the following congruences hold:
\begin{equation}\label{subalgebra condition}\tag{SUB}
 \forall i\in\{1,2,3\}: \; D
\alpha^{-1}\mcR(\alpha[i])(\alpha^{-1})^\tut D\equiv 0\mod D_i.
\end{equation}
By $\alpha[i]$ we denote the $i$-th column of $\alpha\in\G$,
and $(\alpha^{-1})^\tut$ is the transpose of~$\alpha^{-1}$.  The
antisymmetric $3\times3$-matrix of $\Zp$-linear forms
$$\mcR(\bfy):=(\tuL_{ij}({\bfy}))\in\text{Mat}_3(\Zp[\bfy]),$$ where
 $\tuL_{ij}({\bfy}):=\lambda_{ij}^1 y_1+\lambda_{ij}^2
 y_2+\lambda_{ij}^3 y_3$, encodes the structure constants
 $\lambda_{ij}^k$ of the algebra $L$ with respect to the given
 $\Zp$-basis $(e_1,e_2,e_3)$.

 Our opening remarks now amount to observing that
 condition~\eqref{subalgebra condition} is satisfied for all values of
 $r_0$ greater than or equal to the minimal value of $r_0$ with this
 property. The task of determining the summand $|L:\Lambda_0|^{-s}$ of
 $A(s)$ thus reduces to the problem of calculating this minimal value.
 For each subset $I\subseteq\{1,2\}$, we set
$$A_I(s):=\sum_{I([\Lambda])=I}|L:\Lambda_0|^{-s},$$ so that
 $$A(s)=\sum_{I\subseteq\{1,2\}}A_I(s).$$ Our aim is to compute the
Dirichlet series $A_I(s)$ by investigating the
condition~\eqref{subalgebra condition} in each of the four cases.

\medskip \ni \emph{Case} $I=\varnothing$. Clearly in this case the
condition~\eqref{subalgebra condition} is trivially satisfied for all
$r_0\in\N_0$. As there is a unique homothety class of type
$\bfr=(0,0)$, we obtain $A_{\varnothing}(s)=1$.

\medskip
\ni \emph{Case} $I=\{1\}$. In this case
$D=p^{r_0}\textup{diag}(p^{r_{1}},1,1)$, $r_1\in\N$. We have
$$\G_{(r_1,0)}=\left(\begin{array}{cc}\Zp^*
      &\begin{array}{cc}\Zp&\Zp\end{array}\\\begin{array}{c}p^{r_1}\Zp\\p^{r_1}\Zp\end{array}&\GL_2(\Zp)\end{array}\right).$$
      Lattice classes of type $\bfr=(r_1,0)$ may thus be identified with
      the $(p^{-2}+p^{-1}+1)p^{2r_1}$ points of the finite projective
      space $\mathbb{P}^2(\Z/p^{r_1}\Z)\cong\G/\G_{(r_1,0)}$ by taking
      the first column of a matrix modulo units
      in~$\Zp/(p^{r_1}\Zp)$. Only the first of the three conditions
      in~\eqref{subalgebra condition} is non-trivial:
\begin{equation*}
  D\alpha^{-1}\mcR(\alpha[1])(\alpha^{-1})^\tut D\equiv 0\mod
  p^{r_0+r_1}.
\end{equation*}
Since the matrix~$\mcR$ is antisymmetric, we only need to check a
single entry of the matrix on the left hand side. Cancelling a factor
$p^{r_0}$ we obtain that~\eqref{subalgebra condition} holds if and
only if
\begin{equation*}
  p^{r_0}\left(\alpha^{-1}\mcR(\alpha[1])(\alpha^{-1})^\tut\right)_{23}\equiv0\mod p^{r_1}.
\end{equation*}

\begin{lemma}
 For $\mcR(\bfy)=(\tuL_{ij}(\bfy))$ and $\alpha=(\alpha_{ij})$ we have
 $$\det(\alpha)(\alpha^{-1}\mcR(\alpha[1])(\alpha^{-1})^\tut)_{23}=\tuL_{23}(\alpha[1])\alpha_{11}-\tuL_{13}(\alpha[1])\alpha_{21}+\tuL_{12}(\alpha[1])\alpha_{31}.$$
\end{lemma}

\begin{proof}
 Set $\alpha^{-1}=\beta=(\beta_{ij})$. Note that
 $\alpha=\det(\alpha)\beta^\natural$, where $\beta^\natural$ denotes the
 adjoint matrix of~$\beta$. It is now easy to verify that 
 \begin{multline*}
(\alpha^{-1}\mcR(\alpha[1])(\alpha^{-1})^\tut)_{23}=\tuL_{23}(\alpha[1])(\beta_{22}\beta_{33}-\beta_{23}\beta_{32})-\\\tuL_{13}(\alpha[1])(-\beta_{21}\beta_{33}+\beta_{23}\beta_{31})+\tuL_{12}(\alpha[1])(\beta_{21}\beta_{32}-\beta_{22}\beta_{31})\\
=\left(\tuL_{23}(\alpha[1])\alpha_{11}-\tuL_{13}(\alpha[1])\alpha_{21}+\tuL_{12}(\alpha[1])\alpha_{31}\right)/\det(\alpha)
 \end{multline*}
as required.
\end{proof}
We remark that this is a very special situation: for a quadratic
polynomial in the $\beta_{ij}$ to be a linear polynomial in the
entries of the matrix $\alpha\in\GL_d(\Zp)$ it is necessary that
$d=3$.

As $\det(\alpha)$ is a $p$-adic unit, the lemma shows that in the
present case it suffices to control the values taken by the $p$-adic
valuation of the single ternary quadratic form
\begin{equation}\label{form}
f(\bfx):=\tuL_{23}(\bfx)x_1-\tuL_{13}(\bfx)x_2+\tuL_{12}(\bfx)x_3.
\end{equation}
More precisely,
$$\eqref{subalgebra condition}\quad\Leftrightarrow\quad r_0 \geq
r_1-v_p(f(\alpha[1])).$$ Therefore
\begin{equation}\label{A1 formula}
 A_{\{1\}}(s)=\sum_{r_1=1}^\infty\sum_{m=0}^\infty\mcN_{r_1,
 m}\left(p^{-s}\right)^{r_1+3(r_1-m)},
\end{equation}
with
$$\mcN_{r_1,m}:=|\{\bfx=(x_1:x_2:x_3)\in\Proj^2(\Z/p^{r_1}\Z)\mid
\min\{r_{1},v_p(f(\bfx))\}=m\}|. $$ Note that $\mcN_{r_1,m}=0$ unless
$0\leq m\leq r_1$. We shall compute $A_{\{1\}}(s)$ in terms of
$Z_f(s)$ in~\eqref{A_1 Z_f} below.

\begin{remark} It is clear that $f$ depends on the structure constants $\lambda_{ij}^k$ of $L$ with respect to the chosen $\Zp$-basis $(e_1,e_2,e_3)$. In fact, if
\begin{equation*}
 A:=\left(\begin{array}{ccc}\lambda_{23}^1&\lambda_{31}^1&\lambda_{12}^1\\
\lambda_{23}^2&\lambda_{31}^2&\lambda_{12}^2\\
\lambda_{23}^3&\lambda_{31}^3&\lambda_{12}^3\\
\end{array}\right)
 \end{equation*}
then $f(\bfx)=\bfx A \bfx^\tut$. If we change the basis
$(e_1,e_2,e_3)$ to a basis $(e_1',e_2',e_3')$, where
$e_i'=\sum_{i=1}^3p_{ij}e_j$ for a matrix $P=(p_{ij})\in\G$, the
quadratic form $f$ is transformed into $f'(\bfx)=\bfx A'\bfx^\tut$,
where $A'=(\det P)P^{-1} A (P^{-1})^\tut$
(cf.~\cite{TasakiUmehara/92}). We call ternary quadratic forms $f$ and $f'$ equivalent if they are related in this way.
\end{remark}

\medskip
\ni \emph{Case} $I=\{2\}$. This case is much simpler than the previous
one. Note that $D=p^{r_0}\text{diag}(p^{r_{2}},p^{r_{2}},1)$, $r_2\in\N$. We have
$$\G_{(0,r_2)}=\left(\begin{array}{cc}\GL_2(\Zp)
      &\begin{array}{c}\Zp\\
      \Zp\end{array}\\\begin{array}{cc}p^{r_2}\Zp&p^{r_2}\Zp\end{array}&\Zp^*\end{array}\right).$$

Lattice classes of type $\bfr=(0,r_2)$ may thus be identified with the
$(p^{-2}+p^{-1}+1)p^{2r_2}$ points of the finite Grassmannian
$\textup{G}(2,3)(\Z/p^{r_2}\Z)\cong\G/\G_{(0,r_2)}$, determined by the
first two columns of a matrix modulo $p^{r_2}$. One checks immediately
that the subalgebra condition~\eqref{subalgebra condition} is
satisfied for all $r_0\in\N_0$. Thus
$$A_{\{2\}}(s)=\sum_{r_2=1}^\infty(p^{-2}+p^{-1}+1)\left(p^{2-2s}\right)^{r_2}=(p^{-2}+p^{-1}+1)\frac{p^{2-2s}}{1-p^{2-2s}}.$$
\medskip \ni \emph{Case} $I=\{1,2\}$. We shall see that this case
reduces to the case $I=\{1\}$. We have
$D=p^{r_0}\text{diag}(p^{r_{1}+r_{2}},p^{r_{2}},1)$, where
$r_1,r_2\in\N$.  Therefore
$$\G_{(r_1,r_2)}=\left(\begin{array}{ccc}\Zp^*
 &\Zp&\Zp\\p^{r_1}\Zp&\Zp^*&\Zp\\p^{r_1+r_2}\Zp&p^{r_2}\Zp&\Zp^*\end{array}\right).$$
 One verifies without difficulty that, as in the case $I=\{1\}$,
$$\eqref{subalgebra condition}\quad\Leftrightarrow\quad r_0 \geq
r_1-v_p(f(\alpha[1])).$$ In other words, if the lattice class
$[\Lambda]$ is of type $\bfr=(r_{1},r_{2})\in\N^2$ and given by the
left-coset $\alpha\G_\bfr$, the subalgebra condition~\eqref{subalgebra
condition} only depends on the left-coset
$\alpha\G_{(r_{1},0)}$. Evidently each fibre of the natural projection
$$\G/\G_\bfr\twoheadrightarrow\G/\G_{(r_{1},0)}$$ has
cardinality~$(p^{-1}+1)p^{2r_2}$. The computation of $A_{\{1,2\}}(s)$
reduces therefore to the computation of $A_{\{1\}}(s)$. Indeed,
\begin{align*}
 A_{\{1,2\}}(s)&=\sum_{r_2=1}^\infty(p^{-1}+1)\left(p^{2-2s}\right)^{r_2}\sum_{r_1=1}^\infty\sum_{m=0}^\infty\mcN_{r_1,m}\left(p^{-s}\right)^{r_1+3(r_1-m)}\\&=(p^{-1}+1)\frac{p^{2-2s}}{1-p^{2-2s}}A_{\{1\}}(s).
\end{align*}

\medskip
We now develop a formula for the Dirichlet series $A_{\{1\}}(s)$ in
terms of Igusa's local zeta function $Z_f(s)$ associated to the
ternary quadratic form $f$ given in~\eqref{form}. Writing $r$ for
$r_1$ we express the numbers $\mcN_{r,m}$ in our description~\eqref{A1
formula} of $A_{\{1\}}(s)$ in terms of the integers $N_m^\star$,
defined in Section~\ref{igusa}. First we rephrase the $\mcN_{r,m}$ --
counting solutions of equations in finite \emph{projective} spaces --
in terms of the numbers of solutions in corresponding \emph{affine}
cones. We set
\begin{equation}\label{Nrm starred unstarred}
N_{r,m}^\star:=|\{\bfx\in W(r)\mid \min\{r,v_p(f(\bfx))\}=m\}|,
\end{equation}
observe that $N_{r,m}(1-p^{-1})p^r=N_{r,m}^\star$ and that
\begin{align}\label{formula Nrm starred}
 N_{r,m}^\star=\left\{\begin{array}{ll}\mu_m^\star p^{nr}&\text{ if
 }m<r,\\ N_r^\star&\text{ if }m=r. \end{array}\right.
\end{align}

\ni As $n=3$ in our specific situation, this allows us to write

\begin{align*}
 A_{\{1\}}(s)&=\sum_{r=1}^\infty\sum_{m=0}^r\frac{N^\star_{r,m}}{(1-p^{-1})p^r}(p^{-s})^{4r-3m}\\
             &=\frac{1}{1-p^{-1}}\left(\sum_{r=1}^\infty\sum_{m=0}^{r-1}\mu_m^\star
             p^{2r}(p^{-s})^{4r-3m}+\sum_{r=1}^\infty
             p^{-r}N_r^\star(p^{-s})^r\right)\\ &=
             \frac{1}{1-p^{-1}}\left(\sum_{m=0}^\infty\left(\frac{N_m^\star}{p^{3m}}-\frac{N_{m+1}^\star}{p^{3(m+1)}}\right)t^{-3m}\sum_{r=m+1}^\infty(p^2t^4)^r\right.\\&\quad\quad\quad\left.-p^{-3}\sum_{r=1}^\infty(p^2t^4)^r+P^\star(p^2t)-1\right)\\&=\frac{1}{1-p^{-1}}\left(\left(Z^\star(s-2)+p^{-3}\right)\frac{p^2t^4}{1-p^2t^4}-\frac{p^{-1}t^4}{1-p^2t^4}\right.\\
             &\quad\quad\quad\left.+P^\star(p^2t)-1\right).
\end{align*}
Here we used the identities \eqref{Nrm starred unstarred},
\eqref{formula Nrm starred}, \eqref{mustar formula}, \eqref{Pstar
  formula} and \eqref{Zstar formula} and wrote $t=p^{-s}$. Using
equation~\eqref{Pstar Zstar} (in which we replace $t$ by $p^2t$, i.e.
$s$ by $s-2$), we obtain
\begin{equation}\label{A_1 Z_f}
 A_{\{1\}}(s)=\frac{1-p^{-3}}{1-p^{-1}}\left(\frac{p^2t}{1-p^2t}-\frac{Z^\star_f(s-2)p^2t(1-t^3)}{(1-p^{-3})(1-p^2t^4)(1-p^2t)}\right).
\end{equation}
Notice that $Z_f^\star(s-2)=0$ if $L$ is abelian. Using the
identity~\eqref{Z Zstar} and the observation that
\begin{align*}A_{\{1\}}(s)+A_{\{1,2\}}(s)&= A_{\{1\}}(s)\left(1+(p^{-1}+1)\frac{p^2t^2}{1-p^2t^2}\right)\\&= A_{\{1\}}(s)\frac{1+pt^2}{1-p^2t^2}
\end{align*}
it is now immediate that
\begin{align*}
  \zeta_{L}(s)=\frac{1}{1-t^3}\sum_{I\subseteq\{1,2\}}A_I(s)=\zeta_{\Zp^3}(s)-\frac{Z_f(s-2)p^2t}{(1-p^2t^2)(1-p^2t)(1-p^{-1})}.
\end{align*}
This proves Theorem~\ref{theorem} for $i=0$. But passing from $L$ to
$p^iL$ amounts to replacing $f$ by $p^if$. It is clear, however,
that $Z_{p^if}(s)=t^iZ_f(s)$, so that
$Z_{p^if}(s-2)=(p^2t)^iZ_f(s-2)$. The result follows.

\section{Explicit computations}\label{examples}

In this section we show how Theorem~\ref{theorem} gives rise to simple
computations of the zeta functions of $3$-dimensional $\Zp$-Lie
algebras. We shall tacitly assume the notation from
Section~\ref{igusa}.

\subsection{The Heisenberg Lie algebra $\mathfrak{h}(\Zp)$} \label{heisenberg}

The (local) Heisenberg Lie algebra has a presentation
$$\mathfrak{h}(\Zp)=\Zp e_1\oplus \Zp e_2 \oplus \Zp e_3,$$ where
$[e_1,e_2]=e_3$ is the only non-zero relation. We thus have
$$\mcR(\bfy)=\left(\begin{array}{ccc}&y_3&\\-y_3&&\\&&\end{array}\right).$$
The ternary quadratic form equals $f(\bfx)=x_3^2$. It is well-known
that Igusa's local zeta function associated to $f$ equals
$$Z_f(s)=\int_{\Zp}|x^2|_p^s\tud x=\frac{1-p^{-1}}{1-p^{-1-2s}}.$$
Using Theorem~\ref{theorem} it is now easy to confirm the formula for
the local Heisenberg Lie algebra (cf.~\cite{GSS/88}):
\begin{align*}
 \zeta_{\mathfrak{h}(\Zp)}(s)&=\zeta_{\Zp^3}(s)-\frac{p^2t}{(1-p^2t^2)(1-p^2t)(1-p^3t)}\\&=\zeta_p(s)\zeta_p(s-1)\zeta_p(2s-3)\zeta_p(2s-2)\zeta_p(3s-3)^{-1}.
\end{align*}

\subsection{The `simple' Lie algebra $\mathfrak{sl}_2(\Zp)$}\label{sl2}
The Lie algebra $\mathfrak{sl}_2(\Zp)$ has a presentation
$$\mathfrak{sl}_2(\Zp)=\Zp e_1\oplus \Zp e_2 \oplus \Zp e_3,$$
where 
$$ [e_1,e_2]=e_3,\quad [e_1,e_3]=-2e_1, \quad[e_2,e_3]=2e_2.$$ We obtain
$$\mcR(\bfy)=\left(\begin{array}{ccc}&y_3&-2y_1\\-y_3&&2y_2\\2y_1&-2y_2&\end{array}\right).$$
The relevant ternary quadratic form is thus
$f(\bfx)=x_3^2+4x_1x_2$. Note that, for $p>2$, $f$ defines a smooth
conic in projective $2$-space which has $p+1$ points over $\Fp$ and
good reduction modulo~$p$. It follows from Denef's formula for Igusa's
local zeta function in this case (cf., for instance,~\cite[(6) on
p. 1146]{DenefMeuser/91}) that
\begin{equation}\label{sl2 formula}
 Z_f(s-2)=\frac{(1-p^{-1})(1-p^{-1}t)}{(1-pt^2)(1-pt)}.
\end{equation} We obtain from Theorem~\ref{theorem} the known formula
\begin{align*}
 \zeta_{\mathfrak{sl}_2(\Zp)}(s)&=\zeta_{\Zp^3}(s)-\frac{(1-p^{-1}t)p^2t}{(1-pt)(1-p^2t)(1-pt^2)(1-p^2t^2)}\\&=\zeta_p(s)\zeta_p(s-1)\zeta_p(2s-1)\zeta_p(2s-2)\zeta_p(3s-1)^{-1}
\end{align*}
and, more generally, the formulae for $\zeta_{p^i\mathfrak{sl}_2(\Zp)}(s)$ computed in \cite[Theorem 3.1]{duSsl2/00}.

\medskip
The case $p=2$ is nearly as simple. We shall in fact derive a formula
for $Z_f(s-2)$, $f(\bfx)=x_3^2+p^2x_1x_2$, valid for all primes~$p$,
in terms of the function $Z_{\tilde{f}}(s-2)$, where
$\tilde{f}(\bfx)=x_3^2+x_1x_2$. As $4$ is a $p$-adic unit for odd $p$,
the right hand side of~\eqref{sl2 formula} yields a formula for
$Z_{\tilde{f}}(s-2)$, valid for all primes. To compute $Z_f(s)$ we
define $W_3:=\{\bfx\in\Zp^3\mid x_3\in\Zp^*\}$ and write
\begin{equation*}
 Z_f(s)=\int_{W_3}|f(\bfx)|_p^s\tud\mu + \int_{\Zp^3\setminus W_3}|f(\bfx)|_p^s\tud\mu.
\end{equation*}
The first summand equals $\mu(W_3)=1-p^{-1}$. To compute the second summand, we
perform a change of variable $x_3=px_3'$, say, effectuating a change
of measure $\tud\mu=|p|_p\tud\mu'$. Thus
$$\int_{\Zp^3\setminus W_3}|f(\bfx)|_p^s\tud\mu=\int_{\Zp^3}|p^2\tilde{f}(x_1,x_2,x_3')|_p^s\,|p|_p\tud\mu'=p^{-1-2s}Z_{\tilde{f}}(s).$$
Combining these pieces of information we obtain, for $p=2$,
$$Z_f(s-2)=1-2^{-1}+8\cdot2^{-2s}Z_{\tilde{f}}(s-2).$$ Our
Theorem~\ref{theorem} now confirms that, for $p=2$,
\begin{align*}
\zeta_{\mathfrak{sl}_2(\Z_2)}(s)&=\zeta_{\Z_2^3}(s)-\frac{(1-2\cdot2^{-s}+6\cdot2^{-2s})2^{2-s}}{(1-2^{1-s})(1-2^{2-s})(1-2^{1-2s})(1-2^{2-2s})}\\&=\zeta_2(s)\zeta_2(s-1)\zeta_2(2s-2)\zeta_2(2s-1)(1+6\cdot2^{-2s}-8\cdot2^{-3s}).
\end{align*}
This formula was first given in~\cite{duSTaylor/02,White/00}.

\subsection{The `simple' Lie algebra $\mathfrak{sl}_1(\Delta_p)$}\label{D_p}
In~\cite{Klopsch/03}, Klopsch computed the zeta function of
$L=\mathfrak{sl}_1(\Delta_p)$, where $\Delta_p$ is the maximal
$\Zp$-order in a central $\Qp$-division algebra of index~$2$. The Lie algebra $L$ contains elements $\bfi,\bfj,\bfk$ satisfying the relations
\begin{equation}\label{ijk relations}
 [\bfi,\bfj]=\bfk,\quad [\bfi,\bfk]=\rho\bfj,\quad[\bfj,\bfk]=-p\bfi,
\end{equation}
where $\rho\in\{1,2,\dots,p-1\}$ is a non-square modulo~$p$ if $p$ is odd and $\rho=-3$ if $p=2$. 

For $p>2$, the triple $(e_1,e_2,e_3)=(\bfi,\bfj,\bfk)$ forms a
$\Zp$-basis for $L$. We obtain in this case
$$\mcR(\bfy)=\left(\begin{array}{ccc}&y_3&\rho y_2\\-y_3&&-py_1\\-\rho
y_2&py_1&\end{array}\right)$$ and are thus led to study Igusa's local
zeta function associated to
$$f(\bfx)=x_3^2-\rho x_2^2-px_1^2.$$ The easiest way to do this may be
 to compute the Poincar\'e series $P_f^\star(t)$ (cf.~\eqref{Pstar
 formula}) and then to use the identities~\eqref{Pstar Zstar}
 and~\eqref{Z Zstar}. In fact, the series $P_f^\star(t)$ has only two
 non-zero summands: one easily computes $N_0^\star=1$ and
 $N_1^\star=p-1$. There are, however, no solutions $\bfx$ in $W$ of
 $f(\bfx)\equiv 0$ modulo $(p^2)$. Indeed such a solution would
 necessarily require $x_2\equiv x_3\equiv 0$ modulo $(p)$, forcing
 $x_1\equiv0$ modulo $(p)$. Thus $N_m^\star=0$ for
 $m\geq2$. Lemma~\ref{Pstar lemma} yields
\begin{equation*}
1+(p-1)p^{-3}t=P_f^\star(t)=\frac{1-p^{-3}t-tZ_f^\star(s)}{1-t}
\end{equation*}
and, using \eqref{Z Zstar}, we obtain
\begin{equation*}
  Z_f(s-2)=\frac{(1-p^{-1})(1+p^{-1}+t)}{1-pt^2}.
\end{equation*}
Thus
\begin{align*}
 \zeta_{\mathfrak{sl}_1(\Delta_p)}(s)&=\zeta_{\Zp^3}(s)-\frac{(1+p^{-1}+t)p^2t}{(1-pt)(1-pt^2)(1-p^2t^2)}\\
&=\zeta_p(s)\zeta_p(2s-1)\zeta_p(2s-2).
\end{align*}
\medskip

For $p=2$ the triple
$(e_1,e_2,e_3)$, with
$$e_1=2\bfi,\quad e_2=2\bfj,\quad e_3=\bfj+\bfk,$$ forms a $\Zp$-basis
for $L$. Using~\eqref{ijk relations} we derive the following
commutator relations:
$$[e_1,e_2]=-2e_2+4e_3,\quad[e_1,e_3]=-4e_2+2e_3,\quad[e_2,e_3]=-2e_1.$$
This yields
$$\mcR(\bfy)=\left(\begin{array}{ccc}&-2y_2+4y_3&-4y_2+2y_3\\2y_2-4y_3&&-2y_1\\4y_2-2y_3&2y_1&\end{array}\right)$$
and thus
$$f(\bfx)=-2(x_1^2-2(x_2^2-x_2x_3+x_3^2)).$$
We shall in fact derive a recursion formula for $Z_f(s)$, where
$$f(\bfx)=p(x_1^2-pq(x_2,x_3))$$ and $q$ is an arbitrary binary
quadratic form which is anisotropic modulo~$p$. This condition on $q$
and $p$ is certainly satisfied for the form
$q(x_2,x_3)=x_2^2-x_2x_3+x_3^2$ and the prime $p=2$. The calculation
involves the zeta function $Z_{\tilde{f}}(s)$ associated to the form
$$\tilde{f}(\bfx)=px_1^2-q(x_2,x_3)$$ and uses changes of variables
similar to the one performed in Section~\ref{sl2} for $p=2$. Setting
$W_1:=\{\bfx\in\Zp^3\mid x_1\not\equiv0\mod(p)\}$ and
$W_{2,3}:=\{\bfx\in\Zp^3\mid (x_2,x_3)\not\equiv(0,0)\mod(p)\}$ we
obtain
\begin{align*}
t^{-1}Z_f(s)&=\int_{W_1}|p^{-1}f(\bfx)|_p^s\tud\mu+\int_{\Z_p^3\setminus
W_1}|p^{-1}f(\bfx)|_p^s\tud\mu\\ &=1-p^{-1}+p^{-1}tZ_{\tilde{f}}(s)\\
&=1-p^{-1}+p^{-1}t\left(1-p^{-2}+\int_{\Z_p^3\setminus
W_{2,3}}|\tilde{f}(\bfx)|_p^s\tud\mu\right)\\
&=1-p^{-1}+p^{-1}t\left(1-p^{-2}+p^{-2-s}Z_{p^{-1}f}(s)\right)\\
&=1-p^{-1}+p^{-1}t\left(1-p^{-2}+p^{-2}Z_f(s)\right).
\end{align*}
From this we compute
\begin{equation*}%\label{sl1 p=2}
Z_f(s-2)=\frac{(1-p^{-1})(1+(p+1)t)p^2t}{1-pt^2},
\end{equation*}
and thus, setting $p=2$,
\begin{align*}
\zeta_{\mathfrak{sl}_1(\Delta_2)}&=\zeta_{\Z_2^3}(s)-\frac{(1+3\cdot2^{-s})2^{4-2s}}{(1-2^{2-s})(1-2^{1-2s})(1-2^{2-2s})}\\
&=\zeta_2(s)\zeta_2(2s-1)\zeta_2(2s-2)(1+6\cdot2^{-s}+6\cdot2^{-2s}-12\cdot2^{-3s}),
\end{align*}
confirming the results of~\cite[Theorem 1.1]{Klopsch/03}.

\subsection{Soluble Lie algebras}\label{section soluble}

In this section let $p\geq3$. The soluble $3$-dimensional $\Zp$-Lie
algebras have been listed by Gonz\'alez-S\'anchez and Klopsch
in~\cite{KlopschGonzalez-Sanchez/07}, using an analysis of conjugacy
classes in $\text{SL}_2(\Zp)$.  It suffices to consider the following
maximal representatives of the respective homothety classes listed
below. All others are obtained by multiplying one of the matrices of
relations $\mcR(\bfy)$ by a power of $p$; the effect of this operation
on the zeta function~$Z_f(s)$ is easily controlled. In the following
we choose a notation similar to \cite[\S
6]{KlopschGonzalez-Sanchez/07}. The soluble Lie algebras to be considered are the
following.
\begin{enumerate}
\item The abelian Lie algebra $L_0(\infty)$.

\item The Heisenberg Lie algebra $L_0(0)=\mathfrak{h}(\Zp)$.

\item The non-nilpotent Lie algebra $L_1(0)$. We obtain
    $$
    \mathcal{R}(\mathbf{y}) =
    \begin{pmatrix}
          & - y_2 & - y_3 \\
      y_2 &       & \\
      y_3 &       &
    \end{pmatrix}
    $$
    and
    $$
    f(\bfx) =  0.
    $$
  \item The non-nilpotent Lie algebras $L_2(0,r,d)$ with $r \in \N$ and $d
    \in \Z_p$. We obtain
    $$
    \mathcal{R}(\mathbf{y}) =
    \begin{pmatrix}
                     & - y_2 - p^r d y_3 & - p^r y_2 - y_3 \\
      y_2 + p^r d y_3 &       & \\
      p^r y_2 + y_3   &       &
    \end{pmatrix}
    $$
    and
    $$
    f(\bfx) =  p^r
    (x_2^2 - d x_3^2).
    $$
  \item The non-nilpotent Lie algebras $L_3(0,r,d)$ with $r \in \N_0$ and
    $d \in \Z_p$. We obtain
    $$
    \mathcal{R}(\mathbf{y}) =
    \begin{pmatrix}
                   & - d y_3 & - y_2 - p^r y_3 \\
      d y_3        &         & \\
      y_2 + p^r y_3 &         &
    \end{pmatrix}
    $$
    and
    $$
    f(\bfx) =  x_2^2 + p^r x_2
    x_3 - d x_3^2.
    $$
  \item The non-nilpotent Lie algebras $L_4(0,r)$ with $r \in \N_0$. We obtain
    $$
    \mathcal{R}(\mathbf{y}) =
    \begin{pmatrix}
              & - p^r y_3 & - y_2 \\
      p^r y_3 &           & \\
      y_2     &           &
    \end{pmatrix}
    $$
    and
    $$
    f(\bfx) = x_2^2 - p^r x_3^2.
    $$ 
  \item The non-nilpotent Lie algebras $L_5(0,r)$ with $r \in \N_0$. We obtain
    $$
    \mathcal{R}(\mathbf{y}) =
    \begin{pmatrix}
                  & - p^r \rho y_3 & - y_2 \\
      p^r \rho y_3 &               & \\
      y_2          &               &
    \end{pmatrix},
    $$
    where $\rho\in\Z_p^*$ is a non-square modulo $p$, and
    $$
    f(\bfx) = x_2^2 - p^r \rho x_3^2.
    $$ 
\end{enumerate}
The cases A to C have already been treated. Rather than calculate the
zeta function $Z_f(s)$ in each of the remaining cases, we note that,
after (possibly) dividing by a power of $p$, completing the square and
a coordinate change, this reduces to the computation of the zeta
function for the polynomial $f(x_2,x_3)=x_2^2-dx_3^2$, $d\in\Zp$. We
distinguish two cases. If $d=p^k\rho$, where $\rho\in\Zp^*$ is a
non-square modulo $p$, we define $Z_{\boxslash,k}(s):=Z_f(s)$. If
$d=p^ku^2$, where $u\in\Zp^*$, we set $Z_{\boxempty,k}(s):=Z_f(s)$. Both
cases are easily computed using
\begin{gather*}
Z_{\boxslash,0}(s)=\frac{1-p^{-2}}{1-p^{-2}t^2}, \quad
Z_{\boxempty,0}(s)=\left(\frac{1-p^{-1}}{1-p^{-1}t}\right)^2, \\
Z_{\boxslash,1}(s)=Z_{\boxempty,1}(s)=\frac{1-p^{-1}}{1-p^{-1}t}
\end{gather*}
and the fact that both sequences satisfy the same simple recursion
equation of length two. Indeed, for $*=\boxempty$ or $*=\boxslash$, we
have, for $k\in\N_0$,
\begin{equation*}Z_{*,k+2}(s)=p^{-1}t^2Z_{*,k}(s)+1-p^{-1}.
\end{equation*}
An elementary calculation using these observations yields

\begin{proposition}\label{prop} 
Let $L$ be a soluble $3$-dimensional $\Zp$-Lie
algebra associated to one of the families~D to~G. Then, for suitable $k,\iota\in\N_0$ and $*=\boxempty$ or
$*=\boxslash$,
$$\zeta_{L}(s)=\zeta_{\Zp^3}(s)-Z_{*,k}(s-2)\zeta_p(2s-2)\zeta_p(s-2)p^{(2-s)(\iota+1)}(1-p^{-1})^{-1}.$$
The abscissa of convergence $\alpha$ of $\zeta_{L}(s)$ equals $1$ in
all cases. If $*=\boxempty$ and $k$ is even then $\zeta_{L}(s)$ has a
triple pole at $s=1$. In all other cases $\zeta_{L}(s)$ has a double
pole at $s=1$. 
\end{proposition}

\begin{corollary}\label{coro} 
Assume the setting of Proposition~\ref{prop}. For $n\in\N_0$ denote by
$\sigma_n$ the number of subalgebras of $L$ of index at most
$p^n$. Then there are constants $c_1, c_2\in\mathbb{R}$, depending on
$L$, such that, for all $n\in\N_0$
\begin{alignat*}{2}
c_1p^nn^2 &\leq \sigma_n \leq c_2p^nn^2 &&\quad\text{ if
$*=\boxempty$ and $k$ is even,} \\ c_1p^nn &\leq \sigma_n \leq
c_2p^nn &&\quad\text{ otherwise.}
\end{alignat*}
\end{corollary}

%-------------------------------------------------------------------------

\begin{acknowledgements} Much of the research for this article was
  carried out during a visit to the Mathematisches Forschungsinstitut
  Oberwolfach under the Research in Pairs programme from March 18 to
  31, 2007. We thank the institute for its hospitality and support. We
  also acknowledge support from the Nuffield Foundation in form of a
  Newly Appointed Science Lecturer's grant, held by Voll. This paper
  forms part of Voll's Habilitationsschrift.
\end{acknowledgements}

%-------------------------------------------------------------------------
%-------------------------------------------------------------------------

\providecommand{\bysame}{\leavevmode\hbox to3em{\hrulefill}\thinspace}
\providecommand{\MR}{\relax\ifhmode\unskip\space\fi MR }
\providecommand{\MRhref}[2]{%
  \href{http://www.ams.org/mathscinet-getitem?mr=#1}{#2}
}
\providecommand{\href}[2]{#2}

%\bibliographystyle{amsplain}
%\bibliography{thebibliography}

\end{document}